 \documentclass[a4paper, 10pt, conference]{ieeeconf}      

\IEEEoverridecommandlockouts                              

\usepackage{dcolumn}    
\usepackage{eurosym}
\usepackage{graphicx}
\usepackage{subfig}
\usepackage{amsmath}
\usepackage{verbatim}
\usepackage[normalem]{ulem}

\graphicspath{{Images/}}
\newcolumntype{d}[1]{D{.}{.}{#1}}
\usepackage{siunitx} 
\usepackage{fancyhdr,lipsum}
\usepackage{fancyhdr,lipsum}
\fancypagestyle{kbmcsc}{%
 \fancyhf{} 
 
 \fancyhead[L]{SUBMITTED TO IEEE PES GM 2019}
}%

\makeatletter
\let\ps@IEEEtitlepagestyle\ps@kbmcsc
\makeatother

\title{\LARGE \bf Directly Constraining Marginal Prices in Distribution Grids Using Demand-Side Flexibility}

\author{Shantanu Chakraborty$^{1}$, Kyri Baker$^{2}$, Milos Cvetkovic$^{3}$, Remco Verzijlbergh$^{1}$ and Zofia Lukszo$^{1}$ 
\thanks{$^{1}$S. Chakraborty, R.A.Verzijlbergh and Z.Lukszo are with Faculty of Technology, Policy and Management at Delft University of Technology, 2628 BX Delft, The Netherlands
        { s.t.chakraborty at tudelft.nl}}%
\thanks{$^{2}$K.A. Baker is with the Department of Civil, Environmental and Architectural Engineering, University of Colorado at Boulder,
        Boulder, CO 80309, USA
        { kyri.baker at colorado.edu}}%
\thanks{$^{3}$M. Cvetkovic is with the Department of Electrical Engineering, \mbox{Mathematics} and Computer Science at the Delft University of Technology, 2628 CD Delft, The Netherlands
        { m.cvetkovic at tudelft.nl}}%
}

\begin{document}

\maketitle
\thispagestyle{empty}
\pagestyle{empty}

\begin{abstract}
Recently, the volatility associated with marginal prices has increased due to large scale integration of renewable generation. Price volatility is undesirable from a consumer perspective. To address this issue, we present a framework for hedging that uses duality theory for quantifying the amount of demand-side flexibility required for constraining marginal prices to the consumers maximum willingness to pay for electricity. Using our formulation, we investigate the ability of an Energy Storage System (ESS), as a demand-side flexibility source, to hedge against electricity price volatility across a multi-time period horizon while accounting for its inter-temporal constraints. Additionally, we analyze the economical benefit that operating the ESS under information forecasts brings to the consumers.    
\end{abstract}

\section{INTRODUCTION}

The integration of renewable energy sources such as wind and solar poses both advantages and operational challenges. Price volatility is an issue that emerges from the intermittent nature of renewables as well as demand side uncertainty \cite{Astaneh2013}. The increase in price volatility due to renewable generation has been observed around the world such as in Spain \cite{Pereira2017}, Australia \cite{Higgs2015}, Denmark \cite{Rintamaki2017}, Germany \cite{Ketterer2014} and Texas \cite{Woo2011}. While supply-side flexibility options are capable of addressing price volatility issues \cite{Pereira2017},\cite{Rintamaki2017}, they are subject to grid congestion. Alternatively, demand-side flexibility options are well suited for reducing price volatility by complementing the variable generation of renewables. In this paper, we focus on demand-side flexibility that can be leveraged from an Energy Storage System. \par

Previously, regulatory approaches that focus on time-of-use and incentive based demand response \cite{Asadinejad2017} or consumer willingness to shift demand under these programs \cite{Mohajeryami2017} have been investigated to account for price volatility. Additionally, technology based solutions such as the potential of energy storage to reduce price volatility has also been explored \cite{Masoumzadeh2018}. An alternative approach to addressing price volatility is through hedging. Conventionally, hedging in long-term electricity markets has been achieved through Forward Contracts \cite{Kirschen2004}, while recent approaches such as investing in Energy Efficiency have also emerged \cite{Baatz2018}. With relation to hedging in the day-ahead market, a comparative analysis between Forward Contracts, call options and incentivizing consumers for flexibility provision is performed in \cite{Zhou2017}. \par

In this paper, hedging is achieved by constraining marginal prices through the provision of demand-side flexibility. The crux of our approach is that, with the help of duality theory, we can ensure that in a day-ahead market the consumers do not pay higher than their maximum willingness to pay for electricity. This energy flexibility can be availed from an ESS, and the question that we seek to answer is: \textit{What quantity of demand-side flexibility is required from an Energy Storage System to hedge against electricity price \mbox{volatility} over multiple time periods in a distribution grid with \mbox{distributed} generation?} \par

In our previous work \cite{Chakraborty2019} we used the formulation presented in \cite{Baker2016} to quantify the demand-side flexibility required for constraining the marginal prices to the consumers maximum willingness to pay for electricity, while accounting for grid constraints. The main feature of this formulation is that by applying duality theory, explicit constraints can be applied on electricity prices, which is generally an optimization output which is not known a priori. As a result of the addition of these price constraints in the dual formulation, a new variable pertaining to the amount of demand-side flexibility required is introduced in the primal formulation. An additional feature of this formulation is that it simplifies the demand-side bidding process by accepting bids that are based solely on price rather than (\(price,MWh\)) bids. For completeness of this paper, the entire hedging formulation is described in the following sections. \par

However, our previous work assumes that the demand-side flexibility requested could always be provided and the physical constraints of the flexible source were ignored. To address this gap, in this paper we focus on the practicality of our formulation by investigating the optimal operation of the ESS required for hedging against price volatility across multiple time periods in an electricity market. Furthermore, we investigate the additional economical benefit for the consumers that can be generated by considering information forecasts such as prices, generation and demand values in the operation of the ESS. This is achieved through the construct of Model Predictive Control (MPC) which accounts for the physical constraints of the ESS, and the amount of flexibility provided through the charging/discharging of the ESS over multiple time periods. \par

MPC in relation to demand response has previously been used for planning thermal loads under different price schemes to shift loads away from periods of peak demand \cite{Kircher2015}. In \cite{Qureshi2014}, MPC was used to determine the control actions that maximize the financial benefits for a commercial office building that participates in demand response programs. Additionally, through the application of MPC to demand-side flexible resources \cite{Mai2015}, \cite{Zhang2016} services such as the provision of flexible power, ancillary services, and load following capabilities can be enabled for generators. However, the application of MPC for controlling demand-side flexible resources for hedging against price volatility has received limited attention. \par

There are two contributions of this paper, first we demonstrate the efficacy of our formulation to hedge against price volatility while accounting for inter-temporal constraints. Second, we extend this formulation under a MPC scheme and quantify the amount of flexibility required from an ESS to constrain marginal prices at each time step, while respecting the physical constraints of the ESS. While this formulation can be extended to the transmission grid, in this paper we focus on distribution grids and the entities connected to it such as small-scale industrial, commercial or residential consumers. Under future institutional arrangements where dynamic pricing is adopted in distribution grid markets, the proposed hedging mechanism will prove beneficial to these entities by safeguarding them against rapid price fluctuations. \par  

The structure of this paper is as follows: first we present the problem formulation that uses duality theory to constrain the marginal prices. Then we explain its extension to the optimal power flow while considering MPC of the ESS. A simulation framework that comprises of a physical structure and an accompanying institutional arrangement is then presented. We then illustrate how demand-side flexibility can be used for achieving the hedging functionality and finally conclude by providing recommendations for future research. 

\section{Power Flow Formulation Using Energy Flexibility For Hedging Against Price Volatility}\label{section:OPF}

In this Section, we first provide a quick review of the formulation that we use for constraining the marginal prices required for hedging against price volatility. It is assumed that consumers can leverage demand-side flexibility through flexible resources.

\subsection{Price Hedging Formulation Using Duality Theory}\label{subsection:EDF}
Consider an economic dispatch formulation in which a consumer derives an utility of \euro b/MWh from consuming energy and pays a price of \euro a/MWh in return. It is assumed that the consumer as an entity has a desire to maximize its utility function through flexibly consuming energy within it's lower and upper bounds denoted by \(\underline{P_L}\) and \(\overline{P_L}\) respectively. If we are to neglect the line limit constraints, then the formulation is expressed as follows:

\begin{equation}\label{eq:SWEDGeneral}
\begin{aligned}
& \underset{P_G, P_L}{\text{maximize}}
& & bP_L - aP_G \\
& \text{subject to}
& & \underline{P_L} \leq P_L \leq \overline{P_L} \qquad:\mu_1, \mu_2\\
& & & P_G - P_L = 0 \: \enspace \qquad :\lambda\\
& & & P_G, P_L \geq 0
\end{aligned}
\end{equation}

Associated with this formulation are the dual variables of load consumption (\(\mu_1, \mu_2\)) and power balance (\(\lambda\)). Note, since this problem is linear in nature with convex characteristics, Slater's conditions \cite{Boyd2004} will hold, implying that the optimal value of the dual problem is equal to that of the primal problem. In the dual formulation, by constraining the dual variable of marginal price to be lower than or equal to the consumers maximum willingness to pay, expressed as \(\lambda \leq \pi^{des}\), we introduce a new variable in the primal problem. This new variable is defined as the amount of demand-side flexibility that is required for hedging against the electricity price volatility. The modified primal problem is now expressed as follows

\begin{equation}\label{eq:ModifiedSWEDGeneral}
\begin{aligned}
& \underset{P_G, P_L, P^{flexreq}}{\text{maximize}}
& & bP_L - aP_G - \pi^{des}P^{flexreq} \\
& \text{subject to}
& & P_G - P_L + P^{flexreq} = 0 \\
& & & \underline{P_L} \leq P_L \leq \overline{P_L} \\
& & & P_G, P_L, P^{flexreq} \geq 0
\end{aligned}
\end{equation} 

For further information about the formulation for constraining the marginal prices, we refer the reader to \cite{Baker2016}. 

\subsection{Hedging Formulation Considering Network Flows and Energy Storage Systems}

For determining the amount of flexibility \(P^{flex}\) achieveable for constraining marginal prices while accounting for network and ESS constraints, the following formulation is used:

\begin{subequations} \label{eq:MPCOPF}
\begin{equation}
\begin{aligned}
\underset{P_G, P_L, P^{flex}}{\text{max}}
\sum^{h=H}_{h=1}b^{load}(t,h)(P^{load}_L)(t,h) - a^{trans}(t,h)P^{trans}_G(t,h)\\ - a^{dist}(t,h)P^{dist}_G(t,h) - \pi^{des}(t,h)P^{flex}(t,h) \forall t=1,...,24
\end{aligned}\tag{\ref{eq:MPCOPF}}
\end{equation}
subject to:
\begin{equation}
P_{G_i}(t,h) - P_{L_i}(t,h) = \sum_{j \in \Omega_i}\frac{\theta_i(t,h) - \theta_j(t,h)}{X_{ij}}, \quad \forall i \in N \backslash K\label{eq:MPCOPFA} 
\end{equation}
\begin{equation}
\begin{split}
  P_{G_k}(t,h) - P_{L_k}(t,h) + P^{flex}_k(t,h) = \\\sum_{j \in \Omega_k}\frac{\theta_k(t,h) - \theta_j(t,h)}{X_{kj}} \enspace \forall k \in K
  \end{split}\label{eq:MPCOPFB}
\end{equation}
\begin{equation}
    \theta_1 = 0 \label{eq:MPCOPFC}
\end{equation}
\begin{equation}
    \underline{P_{L_i}} \leq P_{L_i}(t,h) \leq \overline{P_{L_i}} \quad \forall i \in N \label{eq:MPCOPFD}
\end{equation}
\begin{equation}
    -\overline{P_{ij}} \leq \frac{\theta_i(t,h) - \theta_j(t,h)}{X_{ij}} \leq \overline{P_{ij}} \quad \forall (i,j) \in \Omega_{ij} \label{eq:MPCOPFE}
\end{equation}
\begin{equation}
    0 \leq P_{G_i}(t,h) \leq \overline{P_{G_i}} \label{eq:MPCOPFF} \quad \forall i \in N
\end{equation}
\begin{equation}
    -P^{flex}_k \leq P^{flex}_k(t,h) \leq P^{flex}_k \quad \forall k \in K \label{eq:MPCOPFG}
\end{equation}
\begin{equation}
    b_k(t,h) = b_k(t,h-1) - b_k^{loss} - P^{flex}_k(t,h) \quad \forall k \in K \label{eq:MPCOPFH}
\end{equation}
\begin{equation}
0 \leq b_k(t,h) \leq b^{max} \quad \forall k \in K
     \label{eq:MPCOPFI}
\end{equation}

\end{subequations}

In Equation \eqref{eq:MPCOPF} the wholesale market prices and the cost of the distribution grid based generation is given by \(a^{trans}\) and \(a^{dist}\) respectively while the utility is given by \(b^{load}\). The set \(N\) is used to represent all the nodes in the system, while set \(K\) is a subset of \(N\) and is used to represent the set of price constrained nodes. Additionally, set \(\Omega_{i}\) and \(\Omega_k\) are used to represent the neighbour of buses \(i\) and \(k\) respectively and \(\Omega_{ij}\) is the set of all the edges in the network. The variable \(P_G = [P_G^{trans}, P_G^{dist}]\) and in Equation \eqref{eq:MPCOPFA} is used to represent the power balance equation at all nodes other than the price constrained nodes. For representing the power balance at the price constrained nodes, we refer to Equation \eqref{eq:MPCOPFB}. Using Equation \eqref{eq:MPCOPFC}, we express the slack bus. Equations \eqref{eq:MPCOPFD}-\eqref{eq:MPCOPFF} express the bounds on load consumption and generation across all the nodes in the set \(N\) and the line flow limits on all the lines in the network. ESS is used for providing the flexibility in order to constrain the marginal prices, and is modeled using the formulation provided in \cite{Maffei2014}. Inclusion of the ESS model introduces constraints Equation \eqref{eq:MPCOPFG}-\eqref{eq:MPCOPFI} on the amount of flexibility available thereby providing a feasibility check on the value of \(P^{flexreq}\) derived in Section \ref{subsection:EDF}. Equation \eqref{eq:MPCOPFG} specifies the bounds on the amount of power that can be charged or discharged within a particular time interval. A positive value of \(P^{flex}\) denotes a discharging operation while a negative value indicates the charging operation of the ESS. Through Equation \eqref{eq:MPCOPFH} we represent the inter-temporal constraints of the operation of the ESS and the state of the ESS is capped at a maximum value specified by Equation \eqref{eq:MPCOPFI}.\par

When there is sufficient amount of flexibility (i.e. \(P^{flex}=P^{flexreq}\)) we will be able to constrain the marginal prices. An approach that we investigate for ensuring that this criteria is satisfied over a given time period is through predictive control based on information forecasts as expressed in \cite{Baker2017}. In order to achieve this predictive control approach, we assume different control horizons. It must be noted that for longer horizons the amount of uncertainty associated with the forecast increases. However, the aspect of uncertainty analysis is beyond the scope of this paper. For each selected horizon we compute the optimal charging and discharging operation of the ESS for the entire horizon, and then implement the first step while discarding all the other steps in the horizon. In Equation \eqref{eq:MPCOPF}, the variable \(t=1...24\) and represents the planning horizon over which the multi-period optimal power flow is applied, while \(h=1...H\) represents the receding horizon length. When \(H=1\), it represents the case where the operation of the ESS is planned without any information forecasts. 

\section{Simulation Framework}\label{section:SF}
\begin{figure}[htp]
\centering
\vspace{-0.5em}
\includegraphics[width = 5cm]{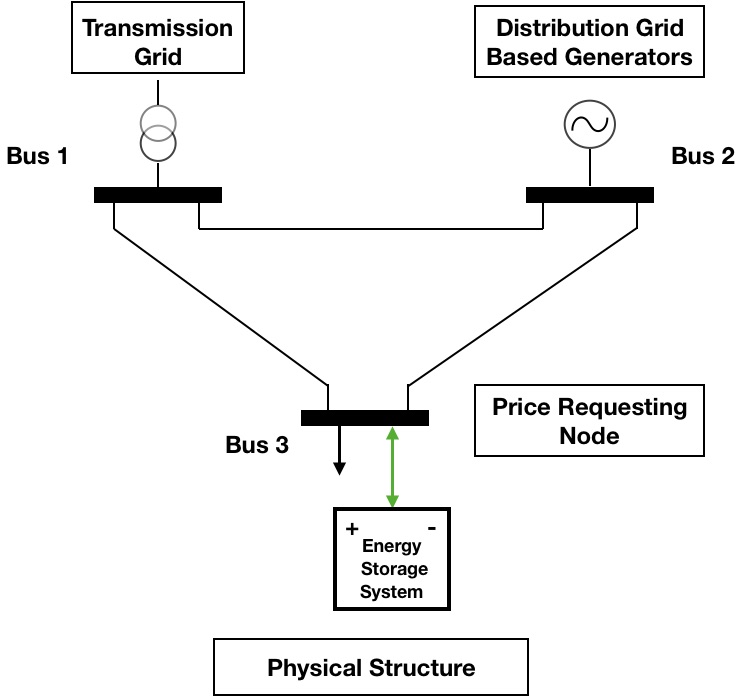}
\vspace{-.2em}
\caption{Simple grid for illustrating proposed hedging strategy}
\label{fig:PS}
\vspace{-1em}
\end{figure}

The simulation framework required for implementing the formulation comprises of two aspects. First, we consider the physical constraints of the network in terms of power flows, generation capacity, energy demand and ESS constraints such as state of charge, losses, and maximum charge/discharge. Then in order to implement the formulation, an institutional arrangement focusing on the coordination of flexibility between different entities at the distribution grid is required. While the formulation presented can be easily extended to larger grids, for illustrative purposes, we assume a simple three bus system as depicted in Figure \ref{fig:PS}. In this figure, the transmission grid is connected to the distribution grid at Bus 1, and given its significantly larger generation capacity compared to the distribution grid based generation, this bus is assumed to be the slack bus. 

Distributed energy resources are connected at Bus 2 in the distribution grid. Consumers who request for a specific price at a bus are known as the Price Requesting Loads and they are connected to the system at Bus 3. Bus 3 in our setup is referred to as the price constrained bus and since it requires demand-side flexibility for constraining the marginal price, an ESS is connected to it. \par

\begin{figure}[htp]
\centering
\vspace{-1em}
\includegraphics[width=6.5cm]{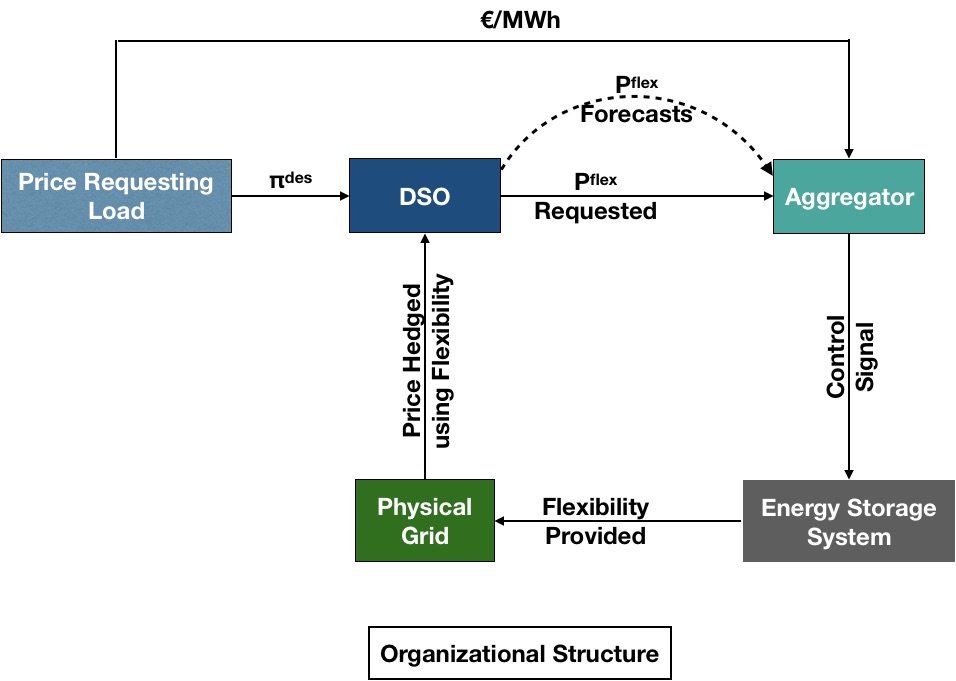}
\vspace{-.5em}
\caption{Organizational Structure for hedging \mbox{strategy}}
\label{fig:OS}
\vspace{-1em}
\end{figure}

Associated with the physical structure of our formulation, we introduce an Organizational Structure, depicted in Figure \ref{fig:OS}, in order to capture actor interactions and information flow between them. The Price Requesting Load communicates its maximum willingness to pay for electricity \(\pi^{des}\) directly to the DSO. This information flow is proposed for two reasons. First, in future local energy markets, the DSO is expected to be a regulated market facilitator at the distribution grid \cite{Gerard2016}. Second, in such an arrangement the DSO is the only entity that has information about cost of generation and physical structure of the network and is thus able to execute the optimal power flow problem. 

By executing the modified optimal power flow formulation that accounts for \(\pi^{des}\) in Equation \eqref{eq:MPCOPF}, the DSO is able to determine the amount of demand-side flexibility required to constrain the marginal prices. When marginal price is lower than \(\pi^{des}\), no additional flexibility is required. However, when marginal prices are higher than \(\pi^{des}\), the DSO on behalf of the Price Requesting Load, makes a request to an aggregator for provisioning the required flexibility. The aggregator, by controlling an ESS which is connected to the price constrained node is able to provide demand-side flexibility required for constraining the marginal prices. It is assumed that the aggregator, in return for its services is remunerated by the Price Requesting Load. \par

It is important to note that the operation of the ESS at one time step affects the energy capacity available at the next. In the absence of information forecasts, the ESS may be operated in a sub-optimal manner over a given time horizon. However, availing information forecasts from the DSO, enables informed decision-making about the operation of the ESS, which leads to increased economical benefits for consumers. This interaction pertaining to information forecasts is captured in Figure \ref{fig:OS}. \par
Additionally, having a larger ESS enables an aggregator to accommodate a wider range of price requests as it is able to provide a higher amount of demand-side flexibility. However, it must be noted that investment in ESS is subject to risk that must be weighed against the price volatility in the electricity market. \par
\begin{figure}[htp]
\centering
\vspace{-1em}
\includegraphics[width=8cm]{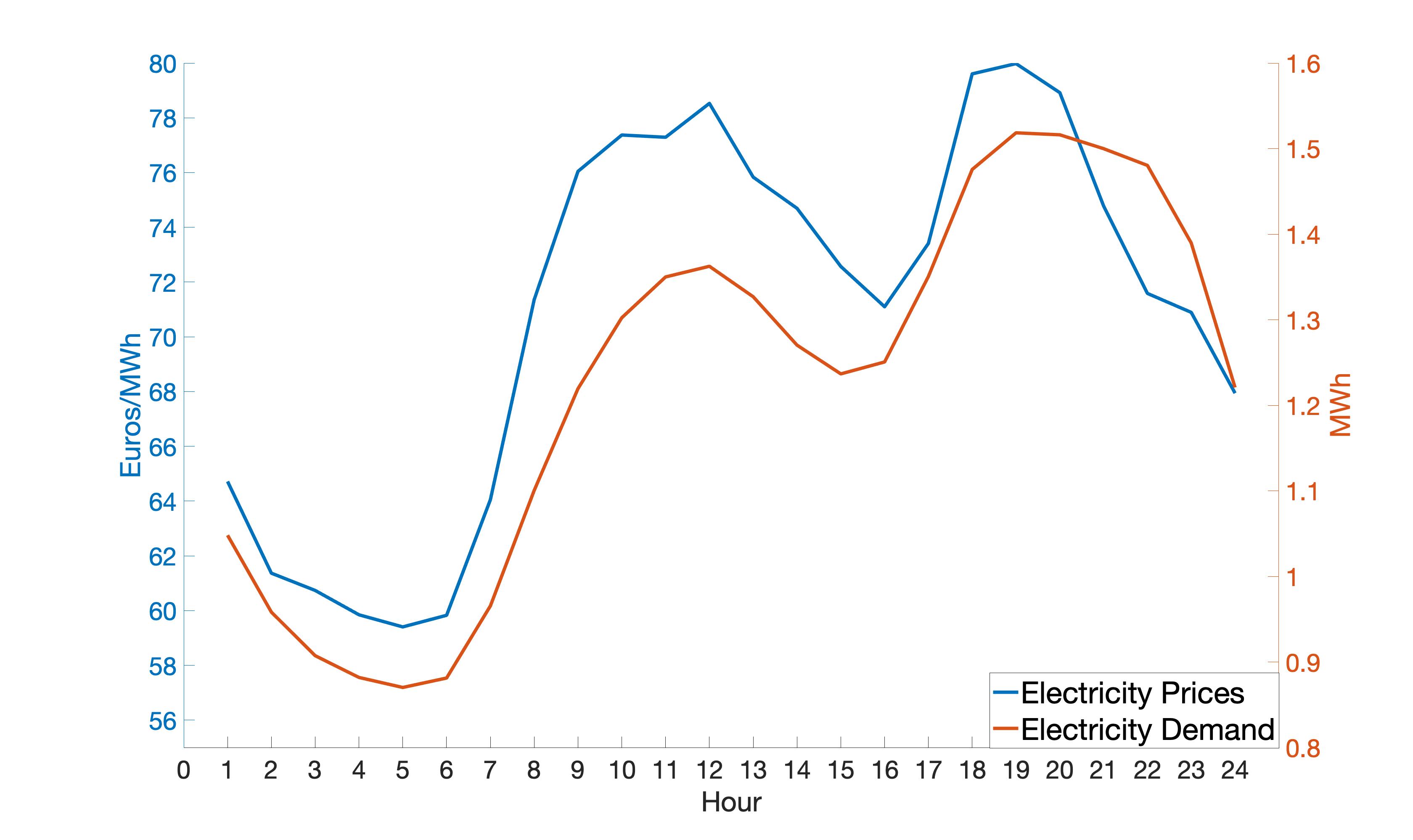}
\vspace{-1em}
\caption{Electricity Prices and Demand}
\label{fig:APXDem}
\vspace{-1em}
\end{figure}

To exhibit the efficacy of our formulation, we focus on the results from a simulation study. We use synthetic data based on the Amsterdam Power Exchange (APX) \cite{ENTSOE}, such that variable \(a^{trans}\) in Equation \eqref{eq:MPCOPF} represents a 24 hour time series. Arbitrary values are selected for the cost of distributed energy sources (\(a^{dist}\)) and load utility values (\(b^{load}\)). The value \(a^{dist}\) is significantly smaller than \(a^{trans}\). The network structure corresponding to Figure \ref{fig:PS} has equal reactances of \(X_{12} = X_{13} = X_{23} = 0.1\Omega\) p.u. Figure \ref{fig:APXDem} provides the Electricity Market Prices and the energy demand used in the simulation.

\section{Results and Discussion}\label{section:SimRes}

\begin{figure}[htp]
\centering
\vspace{-1em}
\includegraphics[width=8cm]{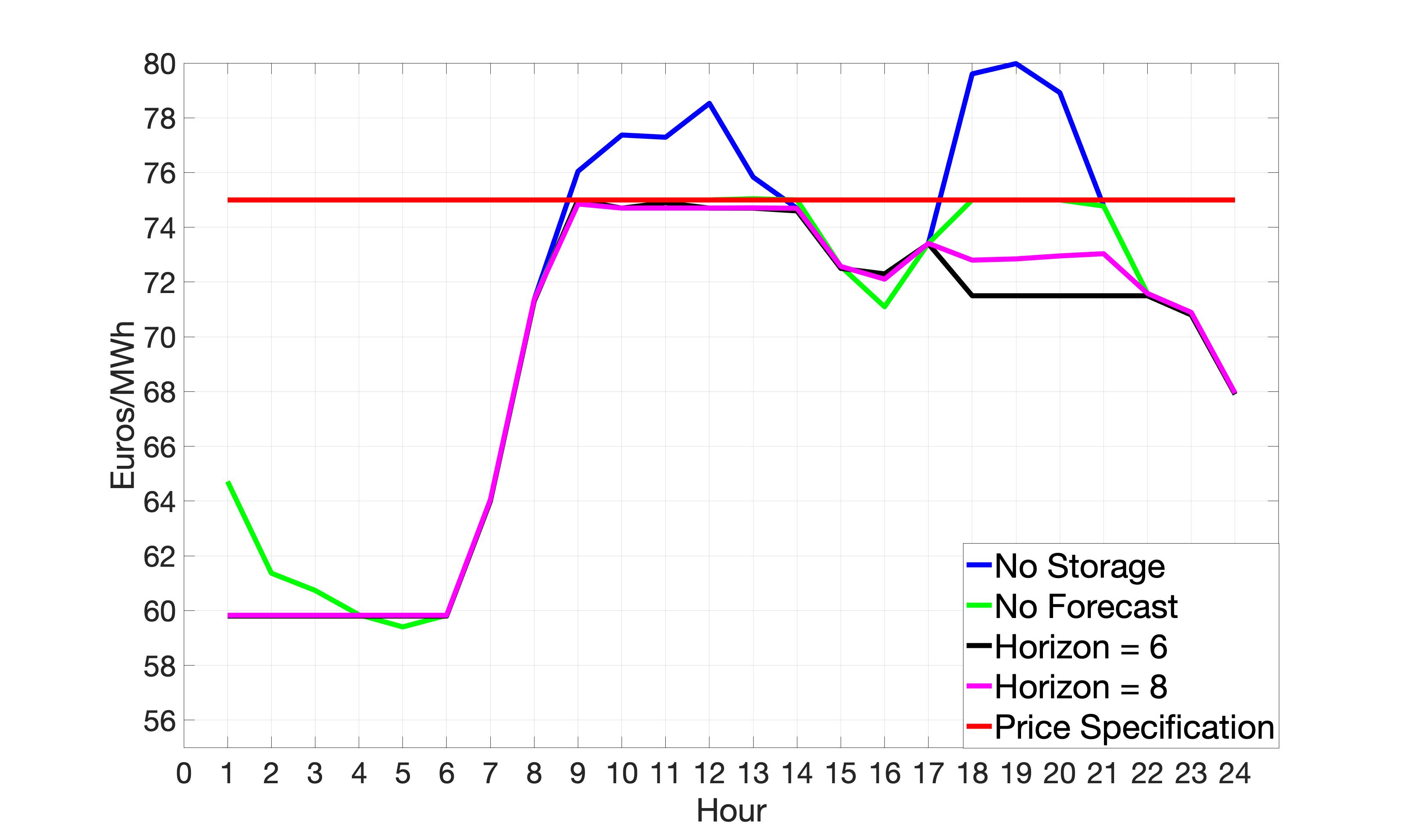}
\vspace{-.5em}
\caption{Constrained LMP under Horizon Lengths 1, 6 and 8}
\label{fig:CompLMP148}
\vspace{-1em}
\end{figure}

In our simulation we assume that the consumers maximum willingness to pay for electricity, \(\pi^{des}=\)\euro 75/MWh. If we compare this maximum willingness to pay against the electricity market prices depicted in Figure \ref{fig:APXDem}, it implies that whenever the prices are expected to rise above \euro 75/MWh, then it should be constrained at this value by availing demand-side flexibility provided from an ESS. Figure \ref{fig:CompLMP148} depicts the effect achieved through our formulation by planning the operation of the ESS under varying forecast horizon lengths of no forecast (green), 6 (black) and 8 (magenta) hours. The consumers maximum willingness to pay, \(\pi^{des}\) (red) is depicted as a constant value over a 24 hour period, and the LMP value without any ESS (blue) is provided as the reference case. Furthermore, in our simulation we consider an ESS of size 2.6MWh and set its initial State of Charge (SoC) to \(75\%\). While the problem of the optimal size of the ESS is out of scope of this paper, as an estimate it should be sized such that it can account for the supply of demand-side flexibility over successive time periods. 

In Figure \ref{fig:CompLMP148}, Hour 9 represents the first instance of when the marginal prices exceeds the value of \(\pi^{des}\), causing the DSO to issue a flexibility request. The marginal prices maintain this status till Hour 13, and during this period demand-side flexibility is provided from the ESS to constrain them. This requirement for flexibility is re-issued by the DSO to the aggregator from Hour 18 - 20 since the LMP again rises above the consumers maximum willingness to pay for electricity.\par

\begin{figure}[htp]
\centering
\vspace{-1em}
\includegraphics[width=9cm]{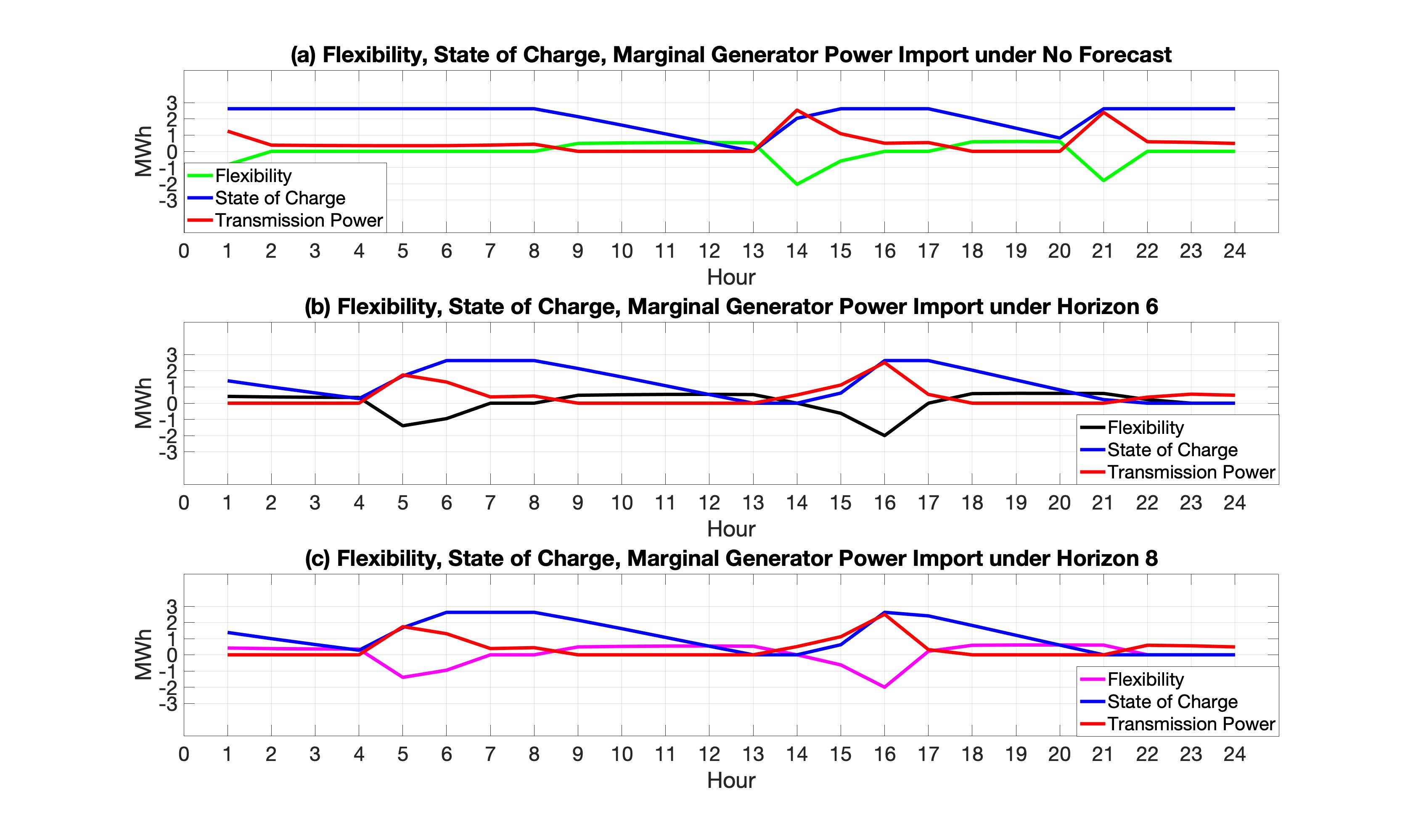}
\vspace{-1em}
\caption{Comparing Flexibility, SoC and Marginal Generation under (a) No forecast, (b) Horizon Length 6 and (c) Horizon Length 8}
\label{fig:FlexSocPg148}
\vspace{-2em}
\end{figure}

In Figure \ref{fig:FlexSocPg148} we visualize  the values of flexibility (charging/discharging of ESS), SoC and the power imported from the transmission grid based generation which is the marginal generator, across the cases when we have no information forecasts Figure \ref{fig:FlexSocPg148}(a), and when the forecast horizons are of length 6 Figure \ref{fig:FlexSocPg148}(b) and 8 Figure \ref{fig:FlexSocPg148}(c) hours. From Hours 1-8, when the LMP values are lower than the specified \(\pi^{des}\), power is imported from the transmission grid. Comparing the MPC horizons against each other in this period, we notice that when horizons are \(H=6\) and \(H=8\) the ESS is discharged during Hours 1-4 which results in lower prices in Figure \ref{fig:CompLMP148} for these horizons as compared to when there is no information forecasts, where the ESS is not discharged. 

At Hour 4, in Figure \ref{fig:FlexSocPg148}, the impact of information forecast is apparent, where we see a portion of the transmission grid based generation being used for pre-charging the ESS when \(H=6\) and \(H=8\). The charging operation of the ESS is inferred from the negative values of flexibility and the subsequent increase in its SoC. From Hours 9-13, the ESS provides flexibility by discharging which results in constraining the marginal prices, providing the hedging functionality. During the Hours 14-16, from Figure \ref{fig:CompLMP148} and Figure \ref{fig:FlexSocPg148}, we can observe the difference that information forecasts can have on charging the ESS. In the absence of information forecasts, the ESS attempts to completely re-charge in the next time-step and hence we notice a relatively higher LMP that is constrained at \(\pi^{des}\). On the other hand, when information forecasts are available, the ESS charging is achieved over time period Hours 14-16 such that the LMP remains lower than \(\pi^{des}\). At Hour 18, the ESS is again discharged for constraining marginal prices till Hour 20. The impact of information forecasts is again observed at Hour 21, when the cost of charging ESS under no forecast is relatively higher. The main difference between MPC horizons \(H=6\) and \(H=8\) are noticed at Hour 17 and Hour 22 with reference to the operation of the ESS resulting in their difference in LMPs as observed in Figure \ref{fig:CompLMP148}. This is an effect of the difference in forecast horizons, as under longer horizons, MPC will take more conservative actions in the current time step. Comparing the economic impact of the information forecasts against no forecast, a saving of \euro 10.24/MWh and \euro 7.76/MWh are observed under horizons \(H=6\) and \(H=8\) respectively over a 24 hour period. 
\section{CONCLUSION}\label{section:Conclusion}
In this paper, a formulation based on duality theory is presented for quantifying the demand-side flexibility required from an Energy Storage System for hedging price volatility in electricity markets. As observed through simulation, flexibility from the ESS is successful in constraining the marginal prices, and over a 24 hour period provides consumers a saving of \euro 23.53/MWh. Under information forecast horizons of 6 and 8 hours, the savings are further increased by \(43.5\%\) and \(33\%\) respectively. Furthermore, the size of an ESS and the length of the forecast horizon have a critical impact on the ability of the ESS to constrain marginal prices and thus the economical benefits provided by them. \par

The results generated from this paper can be extended in multiple directions. First, the optimal power flow method presented here can be modified to account for line losses and increased resistance that are characteristic of low and medium voltage distribution grids. Second, the information forecasts are assumed to be completely deterministic and hence future work will focus on accounting for uncertainties. Finally, the aspect of the optimal size and investment in the ESS was not considered, which could provide improved insights in the future. 

\section*{ACKNOWLEDGMENT}
This work has received funding from the European Union's Horizon 2020 research and innovation programme under the Marie Sklowdowska-Curie grant agreement No. 675318 (\mbox{INCITE})

\bibliographystyle{ieeetr}
\bibliography{BibliographyUpdated.bib}

\end{document}